\documentclass[12pt]{article}
\usepackage{amsmath, amsthm, amsfonts, amssymb}

\usepackage{mathrsfs}
\numberwithin{equation}{section} 

\newtheorem{thm}[equation]{Theorem}
 
\newtheorem{cor}[equation]{Corollary}
\newtheorem{lem}[equation]{Lemma}
\newtheorem{prop}[equation]{Proposition}

\newtheorem{remark}[equation]{Remark}
\newtheorem{defn}[equation]{Definition}

\theoremstyle{definition}

\theoremstyle{remark}

\title{
\vskip-1cm 
\huge 
 Convex ordering for random vectors using 
 predictable representation 
\\
} 

\author 
{
Marc Arnaudon \hspace{-4mm} \setcounter{footnote}{-1} \footnote{a) 
D\'epartement de Math\'ematiques, 
Universit\'e de Poitiers, 
T\'el\'eport 2 - BP 30179, 
86962 Chasseneuil Cedex, 
France. 
marc.arnaudon@math.univ-poitiers.fr}$^{a)}$
 \quad \quad 
Jean-Christophe Breton 
\hspace{-4mm} 
\setcounter{footnote}{-1} \footnote{b) 
D\'epartement de Math\'ematiques, 
Universit\'e de La Rochelle, 
Avenue Michel Cr\'epeau, 
17042 La Rochelle Cedex, 
France. 
jean-christophe.breton@univ-lr.fr 
}$^{b)}$
 \quad \quad 
Nicolas Privault 
\hspace{-4mm} \setcounter{footnote}{-1} \footnote{c) 
Department of Mathematics, 
City University of Hong Kong, 
Tat Chee Avenue, 
Kowloon Tong, 
Hong Kong. 
nprivaul@cityu.edu.hk}$^{c)}$
}

\def\var{{\mathrm{{\rm Var \ \! }}}}

\def\cov{{\mathrm{{\rm Cov \ \! }}}}
\def\Tr{{\mathrm{{\rm Tr \ \! }}}}

\newcommand{\ee}{\mathbb{E}}

\newcommand{\nn}{\mathbb{N}}
\newcommand{\real}{\mathbb{R}}

\newcommand{\retirer}[1]{$ $\newline  }

\def\FF{\mathcal F}

\def\Supp{{\rm Supp \ \!}}

\def\<{\langle}
\def\>{\rangle}

\def\beq{\begin{equation}}
\def\nneq{\end{equation}}
\def\bdef{\begin{defn}}
\def\ndef{\end{defn}}
\def\bthm{\begin{thm}}
\def\nthm{\end{thm}}
\def\bprop{\begin{prop}}
\def\nprop{\end{prop}}
\def\brmk{\begin{remark}}
\def\nrmk{\end{remark}}
\def\bexa{\begin{exa}}
\def\nexa{\end{exa}}
\def\blem{\begin{lem}}
\def\nlem{\end{lem}}
\def\bcor{\begin{cor}}
\def\ncor{\end{cor}}
\def\bexe{\begin{exe}}
\def\nexe{\end{exe}}

\newenvironment{Proof}{\removelastskip\par\medskip \noindent{\em Proof.} \rm}{\penalty-20\null\hfill$\square$\par\medbreak}
\def\bprf{\begin{Proof}}
\def\nprf{\end{Proof}}
\def\bdes{\begin{description}}
\def\ndes{\end{description}}

\allowdisplaybreaks 

\begin{document}
 
\maketitle 

\begin{abstract} 
 We prove convex ordering results 
 for random vectors admitting a predictable representation 
 in terms of a Brownian motion and a non-necessarily independent 
 jump component. 
 Our method uses forward-backward stochastic calculus and 
 extends the results proved in \cite{kmp} in the one-dimensional 
 case. 
 We also study a geometric interpretation 
 of convex ordering for discrete measures 
 in connection with the conditions set on the jump heights 
 and intensities of the considered processes. 
\end{abstract} 
 
\small 

\baselineskip0.42cm 

\vskip8pt 

\noindent 
\textbf{Keywords}: Convex ordering,
 forward-backward stochastic calculus, deviation inequalities, 
 Brownian motion, jump processes. 

\vskip8pt 

\noindent 
\textbf{2000 MR Subject Classification : } 60E15; 60H05,  60G44, 60G55.

\baselineskip0.7cm 

\normalsize 

\maketitle


\baselineskip0.7cm
\section{Introduction} 
 Given two finite measures $\mu$ and $\nu$ on $\real^d$ 
 we say that $\mu$ is convex dominated by $\nu$, and we write 
 $\mu \preceq_{\rm cx} \nu$, if 
\begin{equation} 
\label{cvx} 
 \int_{\real^d} \phi (x) \mu (dx) 
 \leq 
 \int_{\real^d} \phi (x) \nu (dx) 
\end{equation} 
 for all sufficiently integrable 
 convex functions $\phi : \real^d \to \real$. 
 In case $\mu$ and $\nu$ are the respective probability 
 distributions of two random variables $ F $ and $ G $, Relation~\eqref{cvx} 
 is interpreted as the convex concentration inequality 
$$ 
 \ee [ \phi ( F ) ] 
 \leq 
 \ee [ \phi ( G ) ] 
. 
$$ 

\noindent 
 Such concentration inequalities have applications 
 in mathematical finance where they can be interpreted in terms of bounds 
 on option prices on multidimensional underlying assets, see e.g. 
 \cite{bergenthum} and references therein. 

\noindent 
 If $ F $ is more convex concentrated than $ G $ and $G$ is integrable, then 
 $\ee[ F ]=\ee[ G ]$ as follows from taking successively $\phi(x)=x_i$ and $\phi(x)=-x_i$, 
 $i=1,\ldots ,d$, 
 and this amounts to saying that the distributions of $ F $ and $ G $ have 
 same barycenter. 
 On the other hand, applying \eqref{cvx} to the convex 
 function $y \mapsto \phi_x ( y ) = | \langle y - \ee[ F ] , x \rangle |^2$ shows that 
 the matrix $\var  G  - \var  F $ is positive semidefinite, where 
$$
\var  F  = (\cov ( F_i, F_j))_{1\leq i,j \leq d} \quad 
 \mbox{and}  \quad 
 \var  G  = (\cov ( G_i, G_j))_{1\leq i,j \leq d} 
$$ 
 denote the covariance matrices of $ F $ and of $ G $. 
\\ 

\noindent 
 In case $ F $ and $ G $ are Gaussian random vectors with covariance matrices 
 $\Sigma$ and $\tilde{\Sigma}$, these conditions become necessary and sufficient. 
 More precisely, if $\ee[ F ]=\ee[ G ]$ and 
 $\tilde{\Sigma} - \Sigma$ is positive semidefinite then there exists a 
 centered Gaussian random variable $Z$ with covariance 
 $\tilde{\Sigma} - \Sigma$, independent of $ F $ and 
 such that $ G = F +Z$ in distribution, hence 
 $ F $ is more convex concentrated than $ G $ from Jensen's inequality: 
\begin{equation} 
\label{gc} 
 \ee [ \phi ( F ) ] 
 = 
 \ee [ 
 \phi ( 
 \ee [ 
  F +Z
 \mid 
  F  
 ] 
 ) 
 ] 
 \leq 
 \ee [ 
 \ee [ 
 \phi ( 
  F +Z)
 \mid Z 
 ] 
 ] 
 = 
 \ee [ 
 \phi ( 
  F +Z
 ) 
 ] 
 = 
 \ee [ 
 \phi ( 
  G  
 ) 
 ] 
. 
\end{equation} 
\noindent 
 In this paper we aim at obtaining sufficient conditions for the convex ordering 
 of vector-valued random variables, based on their predictable representation as 
 the sum of a diffusion and a jump part. 
 Our main tool of proof consists in the inequality 
$$ 
 \ee [\phi (M(t)+M^*(t))] \leq \ee [\phi (M(s)+M^*(s) )], 
 \quad 
 0 \leq s \leq t 
, 
$$ 
 for all convex functions $\phi : \real^d \to \real$, 
 where $(M(t))_{t\in \real_+}$ and $(M^*(t))_{t\in \real_+}$ 
 are respectively  a forward and a backward $d$-dimensional
 martingale with jumps and continuous parts whose local characteristics 
 satisfy the comparison inequalities assumed in Theorem~\ref{11.-1} below. 
 Such an inequality has been proved in \cite{kmp} for real-valued 
 random variables. 
 We stress however that the arguments of \cite{kmp} are particular to the 
 one-dimensional case and in general they can not be applied to 
 the vector valued setting considered in this paper, for which specific 
 methods have to be developed. 
\\ 

\noindent 
 Note also that by a classical argument, the application 
 of \eqref{cvx} to $\phi (x ) = \exp ( \lambda \|x\| ) $, $\lambda > 0$, 
 entails the deviation bound 
\begin{equation*}
\label{devbdd} 
 P( \Vert  F  \Vert \geq x ) 
 \leq 
 \inf_{\lambda >0} 
 \ee [ e^{\lambda ( \Vert  F  \Vert -x)} {\bf 1}_{\{ \Vert  F  \Vert \geq x\}}] 
 \leq 
 \inf_{\lambda >0} 
 \ee [ e^{\lambda ( \Vert  F  \Vert -x)} ] 
 \leq 
 \inf_{\lambda >0} 
 \ee [ e^{\lambda ( \Vert  G  \Vert -x)} ] 
, 
\end{equation*} 
 $x>0$, hence the deviation probabilities for $ F $ can be estimated via 
 the Laplace transform of $\Vert  G  \Vert$. 
\\ 
 
\noindent 
 We will prove the following type of result. 
 Let $(W(t))_{t \in \real_+}$ and $Z(t)=(Z_1(t), \dots, Z_n(t))$ 
 be respectively a standard $n$-dimensional Brownian motion 
 and a vector of independent real point processes 
 with compensator $(\lambda_1(t), \ldots ,\lambda_n(t))_{t\in \real_+}$ generating a filtration ${\cal F}^M$. 
 Let ${\cal M}^{d\times n}$ denote the set of $d\times n$ 
 real matrices, with ${\cal M}^{d} = {\cal M}^{d\times d}$. 
 Consider $F$ and $G$ two random variables with the 
 predictable representations 
$$
 F = \int_0^\infty A(t) dW(t) +\int_0^\infty J(t) (dZ(t) -\lambda(t) dt)
$$
 where $(A(t))_{t\in \real_+}$, $(J(t))_{t\in \real_+}$ are 
square-integrable ${\cal M}^{d\times n}$-valued $\FF^M_t$-predictable processes, 
 and 
$$
 G= \int_0^\infty \hat{A} (t) d\tilde W(t) +\int_0^\infty \hat{J} (t) (d\tilde Z(t)-\tilde \lambda(t) dt) 
$$
 where $(\hat{A} (t))_{t\in \real_+}$, $(\hat{J} (t))_{t\in\real_+}$ are 
 ${\cal M}^{d\times n}$-valued 
 square-integrable $\FF^M_t$-predictable processes 
 and 
 $\tilde W (t)$ and $\tilde Z(t)=(\tilde Z_1(t), \dots, \tilde Z_n (t))$, 
 $t\in \real_+$, 
 are a $n$-dimensional Brownian motion and a vector of real point processes 
 with respective intensities $(\tilde\lambda_i(t))_{t\in \real_+}$, 
 $i=1,\ldots ,n$, independent of $\FF^M$. 
 In terms of the convex orders cxp, cxpi and psd 
 introduced in Definitions~\ref{order} and \ref{def1} below, 
 we have for example the following corollary of Theorem~\ref{112}. 
 In the sequel, the $\dag$ symbol stands for matrix transposition. 
\bcor 
\label{bt1} 
 The convex concentration inequality \eqref{cvx} holds 
 provided 
$$ 
 A^\dag(t) A(t) \leq_{\rm psd} \hat{A}^\dag(t) \hat{A} (t) 
, 
 \qquad 
 dPdt-a.e., 
$$ 
 and for almost all $t\in \real_+$, we have either: 
\begin{description}
\item{i)}  
$ 
 \displaystyle 
 \sum_{j=1}^n\lambda_j(t) \delta_{(J_{1,j}(t),\dots,J_{d,j}(t))} 
 \preceq_{\rm cxp} 
 \sum_{j=1}^n \tilde\lambda_j(t)  \delta_{(\hat{J}_{1,j}(t), \dots,\hat{J}_{d,j}(t))} 
$
. 
\item or: 
\item{ii)} 
 $J_{i,j}(t) \geq 0$, $\hat{J}_{i,j}(t) \geq 0$, $i=1,\ldots, d$, $j=1,\ldots ,n$, and 
$$ 
 \sum_{j=1}^n\lambda_j(t) \delta_{(J_{1,j}(t),\dots,J_{d,j}(t))} 
 \preceq_{\rm cxpi} 
 \sum_{j=1}^n \tilde\lambda_j(t)  \delta_{(\hat{J}_{1,j}(t), \dots,\hat{J}_{d,j}(t))} 
. 
$$
\end{description} 
\ncor 
\noindent 
 Condition $(ii)$ above will hold in particular if 
$$\lambda_i(t) \leq \tilde{\lambda}_i (t), \qquad 
 \mbox{and} 
 \qquad 
 0\leq J_{i,j}(t) \leq \hat{J}_{i,j}(t),
$$ 
 $1\leq i\leq n$, $1\leq j \leq  d$, for $dt$-almost all $t\in \real_+$. 
 In Theorem~\ref{T5.1}, we provide a geometric interpretation 
 of the convex ordering condition $(i)$ for finitely supported measures. 
\\ 

\noindent 
 In case $F$ and $G$ are Gaussian random vectors with covariance matrices 
 $\Sigma$ and $\tilde{\Sigma}$, we recover \eqref{gc} from Theorem~\ref{bt1} 
 by taking $\lambda (t)=\tilde{\lambda} (t)=0$, 
 $A(t)= {\bf 1}_{[0,T]} (t) \sqrt{\Sigma/T}$ and $\hat{A}(t)= {\bf 1}_{[0,T]} (t) \sqrt{\tilde{\Sigma}/T}$, 
 $t\in \real_+$. 
\\ 

\noindent 
 Note that related convex comparison results 
 have also been obtained in \cite{bergenthum}, \cite{bergenthum2} 
 for diffusions with jumps, under different hypotheses. 
 Namely, it is assumed therein that $G$ is given by the 
 value at time $T$ of a diffusion with jumps. 
 Convex ordering then holds under similar assumptions on the 
 process characteristics, provided the generator 
 of this diffusion satisfies the propagation of convexity property. 
\\ 

\noindent 
 This paper is organized as follows. 
 In Section \ref{s2.0} we introduce the notation 
 of multidimensional forward-backward stochastic calculus 
 with jumps, which will be used in the next sections. 
 In Section~\ref{s2} we prove some convex ordering 
 results for the sums of forward and backward martingales, 
 and in Section~\ref{sc01} we apply those results to 
 random variables given by their predictable representation 
 in terms of a diffusion and a point process. 
 Section~\ref{s5} is devoted to a geometric interpretation 
 of convex ordering for discrete measures on $\real^d$, 
 which gives a better understanding of the conditions set 
 of the jump heights and intensities of the considered point 
 processes. 
\section{Notation} 
\label{s2.0} 
\noindent 
 Let $(\Omega , {\cal F}, P)$ be a probability space equipped with an 
 increasing filtration $(\FF_t)_{t\in \real_+}$ and a decreasing filtration 
 $(\FF^*_t)_{t\in \real_+}$. 
 Consider 
 $(M(t))_{t\in \real_+} 
 = 
 (M_1(t), \dots, M_d(t))_{t\in \real_+}$ 
 a $d$-dimensional $\FF_t$-forward martingale 
 and 
 $(M^*(t) )_{t\in \real_+ } 
 =(M^*_1 (t), \dots, M^*_d (t))_{t\in \real_+ }$ 
 a $d$-dimensional $\FF^*_t$-backward martingale, 
 such that 
 $(M(t))_{t\in \real_+}$ has right-continuous paths with left limits, 
 and 
 $(M^*(t))_{t\in \real_+}$ has left-continuous paths with right limits. 
 Denote respectively by $(M^c(t) )_{t\in \real_+ }$ and 
 $(M^{*c}(t))_{t\in \real_+}$ the continuous parts of 
 $(M(t))_{t\in \real_+}$ and of $(M^*(t) )_{t\in \real_+ }$, and by 
$$\Delta M(t) = M(t) - M(t^-) 
, 
\qquad 
\qquad \Delta^* M^*(t) = M^*(t) - M^*(t^+) 
, 
$$ 
 their forward and backward jumps. 
 The processes $(M(t))_{t\in \real_+}$ and $(M^*(t) )_{t\in \real_+ }$ have jump 
 measures 
$$
\mu (dt ,dx) = \sum_{s>0} {\bf 1}_{\{ \Delta M(s) \not= 0\}} \delta_{(s,\Delta M(s))} (dt ,dx), 
$$ 
 and 
$$
\mu^* (dt ,dx) = \sum_{s>0} {\bf 1}_{\{ \Delta^* M^*(s) \not= 0\}} \delta_{(s,\Delta^* M^*(s))} 
 (dt ,dx), 
$$ 
 where $\delta_{(s,x)}$ denotes the Dirac measure at $(s,x)\in \real_+\times \real^d$. 
 Denote by $\nu (dt , dx)$ and $\nu^* (dt , dx)$ the $(\FF_t)_{t\in \real_+}$ 
 and $(\FF^*_t)_{t\in \real_+}$-dual predictable projections of 
 $\mu (dt ,dx)$ and of $\mu^* (dt ,dx)$. 
 The quadratic variations $([M, M]_t )_{t\in \real_+}$, $([M^* ,M^*]_t )_{t\in \real_+}$ 
 are the ${\cal M}^d$-valued processes defined as the limits in uniform convergence in probability 
$$
 [M,M]_t = \lim_{n\to \infty} 
 \sum_{k=1}^n 
 (M ( t^n_k ) - M ( t^n_{k-1} ) ) (M ( t^n_k ) - M ( t^n_{k-1} ) )^\dag, 
$$ 
 and 
$$
 [M^* ,M^*]_t  = \lim_{n\to \infty} \sum_{k=0}^{n-1} 
 (M^* (t^n_k)-M^* ( t^n_{k+1} ) )(M^* (t^n_k ) -M^* (t^n_{k+1}) )^\dag, 
$$ 
 for all refining sequences 
 $\{0=t_0^n\leq t_1^n \leq \cdots \leq t_n^n=t\}$, 
 $n\geq 1$, of partitions of $[0,t]$ 
 tending to the identity. 
 We let $M^J(t) = M(t) - M^c(t)$, $M^{*J}(t) = M^*(t) - M^{*c}(t)$, 
$$ [ M^J,M^J ]_t = \sum_{0<s\leq t} \Delta M(s)\Delta M(s)^\dag, 
 \qquad 
 [M^{*J},M^{*J}]_t  = \sum_{0\leq s < t} \Delta^* M^* (s) (\Delta^*M^*(s))^\dag, 
$$ 
 with 
$$
\langle M^c,M^c\rangle_t  = [M,M]_t - [M^J,M^J]_t, 
$$ 
 and 
$$
\langle M^{*c},M^{*c} \rangle_t = [M^* ,M^* ]_t - [M^{*J},M^{*J} ]_t, 
$$ 
 $t\in \real_+$. 
 Denote by 
 $(\langle M^J,M^J\rangle_t)_{t\in \real_+}$, $(\langle M^{*J},M^{*J} \rangle_t)_{t\in \real_+}$ 
 the conditional quadratic variations of $( M^J(t))_{t\in \real_+}$ and of $( M^{*J}(t))_{t\in \real_+}$, 
 with 
$$
d\langle M^J,M^J \rangle_t = \int_{\real^d} xx^\dag \nu (dt , dx) 
 \quad 
 \mbox{and} 
 \quad 
 d\langle M^{*J} ,M^{*J} \rangle_t = \int_{\real^d} x x^\dag \nu^* (dt,dx). 
$$ 
 The conditional quadratic variations 
 $(\langle M,M\rangle_t)_{t\in \real_+}$, 
 $(\langle M^*,M^*\rangle_t)_{t\in \real_+}$ 
 of $(M(t))_{t\in \real_+}$ and of $(M^*(t))_{t\in \real_+}$ 
 satisfy 
$$\langle M,M\rangle_t = \langle M^c ,M^c \rangle_t + \langle M^J ,M^J \rangle_t, 
\quad 
 \mbox{and} 
 \quad 
\langle M^*,M^*\rangle_t = \langle M^{*c} ,M^{*c} \rangle_t + \langle M^{*J} ,M^{*J} \rangle_t, 
$$ 
 $t\in \real_+$. 
\noindent 
 Theorem~\ref{11.-1} below and its corollaries are 
 based on the following forward-backward It\^o type change of variable formula, 
 for $(M(t),M^*(t))_{t\in \real_+}$, in which 
 conditions \eqref{hp1} and \eqref{hp2} below are assumed in order to 
 make sense of the integrals with respect to 
 $d M(t)$ and $d^*M^* (t)$. 
 This formula is proved in dimension one in Theorem~8.1 of \cite{kmp}, 
 and its extension to dimension $d\geq 2$ can be proved similarly. 
 For all $f\in {\cal C}^2 (\real^d\times \real^d)$ 
 of the form $f(x,y)=f(x_1, \dots, x_d, y_1, \dots, y_d)$ 
 we have 
\begin{eqnarray} 
\label{itf} 
\lefteqn{ 
 f ( M(t) , M^*(t) )
 = f ( M(s) , M^*(s) )
} 
\\
\nonumber 
& & 
\! \! \! \! \! \! \! \! \! \! 
 + \sum_{i=1}^d\int_{s^+}^t \frac{\partial f}{\partial x_i} (M(u^-) , M^*(u) ) d M_i(u) 
 + 
 \frac12 \sum_{i,j=1}^d\int_s^t 
 \frac{\partial^2 f}{\partial x_i \partial x_j} 
 (M(u) , M^*(u) ) d\langle M^c_i,M^c_j \rangle_u
\\ 
\nonumber 
 & & 
 \! \! \! \! \! \! \! \! \! \! 
 +
 \sum_{ s < u \leq t} 
 \left( 
 f (M(u) , M^*(u) )
 -
 f (M(u^-) , M^*(u) )
 - \sum_{i=1}^d \Delta M_i(u) 
 \frac{\partial f}{\partial x_i} 
 (M(u^-) , M^*(u) )
 \right) 
\\ 
\nonumber 
 & & 
 \! \! \! \! \! \! \! \! \! \! 
 - \sum_{i=1}^d\int_s^{t^-} 
 \frac{\partial f}{\partial y_i} 
 (M(u) , M^*(u^+) )d^*M^*_i(u) 
 - 
 \frac12 
 \sum_{i,j=1}^d \int_s^t 
 \frac{\partial^2 f}{\partial y_i\partial y_j} 
 (M(u) , M^*(u) ) d\langle M^{*c}_i,M^{*c}_j\rangle_u
\\ 
\nonumber 
& & 
 \! \! \! \! \! \! \! \! \! \! 
 - \sum_{ s \leq u < t} 
 \left( 
 f (M(u) , M^*(u) ) - f (M(u) , M^*(u^+) ) 
 - 
 \sum_{i=1}^d\Delta^* M^*_i(u) 
 \frac{\partial f}{\partial y_i} 
 (M(u) , M^*(u^+) ) 
 \right) 
, 
\end{eqnarray} 
 $0\leq s \leq t $, where $d$ and $d^*$ denote the forward and backward It\^o differential, 
 respectively defined as the limits of the Riemann sums 
$$ 
 \sum_{k=1}^n 
 (M_i(t_k^n)-M_i(t_{k-1}^n)) 
 \frac{\partial f}{\partial x_i} 
 (M(t_{k-1}^n) , M^* ( t_{k-1}^n ) ) 
$$ 
 and 
$$ 
 \sum_{k=0}^{n-1} (M^*_i (t_k^n) - M^*_i (t_{k+1}^n) ) 
 \frac{\partial f}{\partial y_i} 
 (M (t_{k+1}^n) , M^* (t_{k+1}^n )) 
$$ 
 for all refining sequences 
 $\{ s =t_0^n\leq t_1^n \leq \cdots \leq t_n^n=t\}$, 
 $n\geq 1$, of partitions of $[s,t]$ tending to the identity. 
\\ 
 
\noindent  
 Here, $\int_0^t \eta (u) dM_i(u)$, resp. $\int_{t}^\infty \eta^*(u) dM_i(u)$, 
 refer to the right, resp. left, continuous version 
 of the indefinite stochastic integrals 
 of the forward, resp. backward, adapted and sufficiently 
 integrable processes $(\eta (u))_{u\in \real_+}$, resp. 
 $ ( \eta^* (u) )_{u\in \real_+}$. 
\section{Convex ordering for martingales} 
\label{s2} 
\noindent 
 We denote by 
 $\langle \cdot,\cdot\rangle$ and $\|\cdot\|$ 
 the usual Euclidean scalar product and norm on $\real^d$. 
 Let ${\cal M}^d$ be the space of real matrices with the scalar product 
$$ 
 \langle A , B \rangle
 : 
 = \Tr(AB^\dag) 
 = \sum_{i,j=1}^d 
 A_{i,j}B_{i,j}, 
 \qquad 
 A,B\in {\cal M}^d 
, 
$$ 
 where we recall that $A^\dag$ stands for the transpose 
 $(A_{j,i})_{1\leq i,j\leq d}$ of $A=(A_{i,j})_{1\leq i,j\leq d}$, 
 and let ${\cal M}_+^d$ be the subset of 
 ${\cal M}^d$ made of positive semidefinite matrices. 
\blem
\label{l31} 
 Let $A$ be a symmetric $d \times d$ matrix. 
 Then the following statements are equivalent: 
\begin{description} 
\item{i)} 
 $A$ is positive semidefinite, 
\item{ii)} 
 for all positive semidefinite matrices $B$ we have 
$\langle A , B \rangle
 \geq 0$, 
\end{description} 
\nlem
\bprf
 Since $A$ is symmetric, if it is positive semidefinite then 
 its spectral decomposition is given as 
$$ 
 A = \sum_{k=1}^d \lambda_k e_k e_k^\dag , 
$$ 
 where the eigenvalues $(\lambda_k)_{k=1,\ldots ,d}$ of $A$ 
 are non-negative and $(e_1, \ldots, e_d)$ denote the eigenvectors 
 of $A$. Hence we have 
$$ 
 \Tr (AB^\dag) 
 = 
 \sum_{k=1}^d 
 \lambda_k  
 \langle e_k , B e_k \rangle 
 \geq 0 
$$ 
 if $B$ is positive semidefinite. 
 The converse follows by choosing $B=x^\dag x$, $x\in \real^d$, and noting that 
$$ 
 \langle x , A x \rangle = \Tr(AB^\dag) \geq 0. 
$$ 
\nprf
\noindent 
\bdef 
\label{def1} 
 Given $A,B\in {\cal M}^d$, we will write $A\leq_{\rm psd} B$ if 
 $B-A$ is positive semidefinite, i.e. 
$$ 
 \langle x , A x\rangle 
 \leq 
 \langle x , B x\rangle, 
 \qquad 
 x\in\real^d. 
$$ 
\ndef 
\noindent
 In the sequel, a function $f:\real^d\to\real$ will be 
 said to be {\em non-decreasing} if 
$$ 
 f(x_1,\ldots,x_d)\leq 
 f(y_1,\ldots,y_d) 
$$ 
 for all $x_1,\ldots,x_d\in\real$ and $y_1,\ldots,y_d\in \real$ 
 such that $x_i\leq y_i$, $i=1,\ldots ,d$. 
 In the sequel, we will need the following orders 
 between positive measures $\mu$, $\nu$ on $\real^d$. 
\bdef 
\label{order} 
~

\begin{description} 
\item{i)} We say that $\mu \preceq_{\rm cxp} \nu$ 
 if 
$$
\int_{\real^d} \phi (x) \mu (dx)\leq\int_{\real^d} \phi (x) \nu (dx)
$$
 for all {\em non-negative} convex functions $\phi :\real^d\to\real^+$. 
\item{ii)} We say that $\mu \preceq_{\rm cxpi} \nu$ if 
$$
\int_{\real^d} \phi (x) \mu (dx)\leq\int_{\real^d} \phi (x) \nu (dx)
$$
 for all {\em non-negative} and {\em non-decreasing} convex 
 functions $\phi :\real^d\to\real^+$. 
\end{description} 
\ndef
\noindent 
 If $\mu$ and $\nu$ are finite measures on 
 $\real^d$, then both $\mu \preceq_{\rm cxp} \nu$ 
 and $\mu \preceq_{\rm cxpi} \nu$ imply 
 $\mu (\real^d) \leq \nu (\real^d)$. 
 More precisely we have the following result. 
\begin{prop} 
 Assume that $\mu$ and $\nu$ are finite measures on $\real^d$. 
 Then $\mu \preceq_{\rm cx} \nu$ is equivalent to 
 $\mu \preceq_{\rm cxp} \nu$ and $\mu (\real^d ) = \nu (\real^d )$. 
\end{prop} 
\begin{Proof} 
 Assume that $\mu \preceq_{\rm cxp} \nu$ and $\mu (\real^d ) = \nu (\real^d )$, 
 and let $\phi \in L^1 (\mu )\bigcap L^1 (\nu )$. 
 For all $a\in \real$ we have 
\begin{eqnarray*} 
\lefteqn{ 
 \int_{\real^d} \phi (x) \mu (dx) 
 - 
 \int_{\{ \phi < a \}} 
 \phi (x) \mu (dx) 
 + a \mu ( \{ \phi < a \} ) 
} 
\\ 
 & = & 
 \int_{\{ \phi \geq a \}} (\phi (x) - a )^+ \mu (dx) 
 + a \mu ( \real^d ) 
\\ 
 & \leq & 
 \int_{\{ \phi \geq a \}} (\phi (x) - a )^+ \nu (dx) 
 + a \nu ( \real^d ) 
\\ 
 & = & 
 \int_{\real^d} \phi (x) \nu (dx) 
 - 
 \int_{\{ \phi < a \}} 
 \phi (x) \nu (dx) 
 + a \nu ( \{ \phi < a \} ) 
, 
\end{eqnarray*} 
 and for $a\leq 0$,
$$ 
 \int_{\{ \phi < a \}} 
 \phi (x) \mu (dx) 
 \leq 
 a \mu 
 ( \{ \phi < a \} 
 ) 
 \leq 0 
, 
\qquad 
 \int_{\{ \phi < a \}} 
 \phi (x) \nu (dx) 
 \leq 
 a \nu 
 ( \{ \phi < a \} 
 ) 
 \leq 0 
, 
$$ 
 hence letting $a$ tend to $-\infty$ yields 
$$ 
 \int_{\real^d} \phi (x) \mu (dx) 
 \leq 
 \int_{\real^d} \phi (x) \nu (dx) 
. 
$$ 
 Conversely we note that 
 $\mu \preceq_{\rm cxp} \nu$ clearly implies $\mu \preceq_{\rm cx} \nu$, 
 and we recover the identity $\mu(\real^d)=\nu(\real^d)$ by 
 applying the property $\mu \preceq_{\rm cx} \nu$ 
 successively with $\phi=1$ and $\phi=-1$. 
\end{Proof} 
\noindent 
 Consequently, $\mu \preceq_{\rm cxp} \nu$ implies 
$$
 \int_{\real^d} 
 x_i 
 \mu (dx) 
 \leq 
 \int_{\real^d} 
 x_i 
 \nu (dx) 
, 
 \quad 
 \mbox{and} 
 \quad 
 - 
 \int_{\real^d} 
 x_i 
 \mu (dx) 
 \leq 
 - 
 \int_{\real^d} 
 x_i 
 \nu (dx) 
, 
 \quad 
 i=1,\ldots , d, 
$$ 
 i.e. $\mu$ and $\nu$ have same barycenter, 
 provided $\mu (\real^d) = \nu (\real^d)$ and $\mu$, $\nu$ 
 are integrable. 
 This also holds when $\mu \preceq_{\rm cxpi} \nu$, 
 $\mu (\real^d) = \nu (\real^d)$ and 
 $\mu$, $\nu$ are supported by $\real_+^d$. 
\\  

\noindent 
 Let now  
\begin{equation} 
\label{hp1} 
 (M(t))_{t\in \real_+}\mbox{ be an }\FF^*_t\mbox{-adapted, }\FF_t\mbox{-forward martingale,} 
\end{equation} 
 and 
\begin{equation} 
\label{hp2} 
 (M^*(t))_{t\in \real_+}\mbox{ be an }\FF_t\mbox{-adapted, }\FF_t^*\mbox{-backward martingale,} 
\end{equation} 
 with characteristics of the form 
\begin{equation*} 
\label{hj} 
\nu (dt,dx) = \nu_t(dx) dt \qquad \mbox{and} \qquad 
 \nu^* (dt,dx) = \nu^*_t (dx) dt 
, 
\end{equation*} 
 and 
\begin{equation*} 
\label{car1} 
d\langle  M^c ,  M^c \rangle_t = H (t) dt, 
 \quad  \mbox{and}  \quad 
d\langle  M^{*c} ,  M^{*c} \rangle_t = H^* (t) dt, 
\end{equation*} 
 where 
 $H(t)=(H_{i,j}(t))_{1\leq i,j\leq d}$ 
 and 
 $H^*(t)=(H^*_{i,j}(t))_{1\leq i,j\leq d}$ 
 are ${\cal M}^d$-valued, $t\in \real_+$, and 
 predictable respectively with respect to 
 $(\FF_t)_{t\in \real_+ }$ and to $(\FF^*_t)_{t\in \real_+ }$. 
 In the sequel, we will also assume that 
 $(H(t))_{t\in \real_+}$, $(H^*(t))_{t\in \real_+} \in L^2 (\Omega \times \real_+)$, 
 and that 
\begin{equation} 
\label{plio} 
 \ee \left[ 
 \int_{\real^d\times \real_+} 
 \Vert x \Vert \nu_t (dx) dt 
 \right] 
 < \infty, 
 \qquad 
 \ee \left[ 
 \int_{\real^d \times \real_+} 
 \Vert x \Vert 
 \nu^*_t (dx) dt 
 \right] 
 < \infty 
. 
\end{equation} 
 The hypotheses on $(H(t))_{t\in \real_+}$ 
 and $(H^*(t))_{t\in \real_+}$ imply that $M^c_t$ 
 and $M^{*c}_t$ are in $L^2(\Omega )$, $t\in \real_+$, 
 and Condition~\ref{plio} is a technical integrability 
 assumption. 
\bthm 
\label{11.-1} 
 Assume that 
$$ 
 H(t) \leq_{\rm psd} H^*(t),  \quad dP dt-a.e.
$$  
 and that for almost all $t\in \real_+$ we have either: 
\begin{description}
\item{i)} 
$\nu_t \preceq_{\rm cxp} \nu^*_t$, 
\item or: 
\item{ii)} 
$\nu_t \preceq_{\rm cxpi} \nu^*_t$ 
 and $\nu_t$, $\nu_t^*$ are supported by $(\real_+)^d$. 
\end{description} 
 Then we have 
\begin{equation} 
\label{holds} 
 \ee [\phi (M(s)+M^*(s) )] 
 \geq 
 \ee [\phi (M(t)+M^*(t))] 
, \qquad 0\leq s \leq t 
, 
\end{equation} 
 for all convex functions $\phi : \real^d \to \real$. 
\nthm 
\noindent 
\noindent 
\begin{Proof} 
 We start by assuming that $\phi$ is a ${\cal C}^2$, convex 
 Lipschitz function and we apply It\^o's formula \eqref{itf} 
 for forward-backward martingales to $f(x,y) = \phi (x+y)$. 
 Taking expectations on both sides of It\^o's formula we get 
\begin{align} 
\nonumber 
& \ee [ \phi ( M(t) + M^*(t) ) ]\\
\nonumber
& = 
 \ee [ \phi ( M(s) + M^*(s) ) ]
 + \frac12 \sum_{i,j=1}^d
 \ee \left[ 
 \int_s^t \phi_{i,j}'' (M(u)+M^*(u) ) 
 d(\langle M^c_i,M^c_j\rangle_u  - \langle M^{*c}_i , M^{*c}_j \rangle_u) 
 \right] 
\\
\nonumber
& +
 \ee \left[
 \int_s^t
 \int_{\real^d}
 ( 
 \phi (M(u) + M^*(u) + x )
 - 
 \phi (M(u) + M^*(u) )
 - \langle 
 x , \nabla \phi (M(u)+M^*(u) ) 
 \rangle 
 ) 
 \nu_u ( dx ) du 
\right] 
\\ 
\nonumber
&  - 
 \ee \left[
 \int_s^t
  \int_{\real^d}
 (  \phi (M(u) + M^*(u)  + x ) -  \phi (M(u) + M^*(u) ) - 
 \langle 
 x , \nabla \phi (M(u)+M^*(u) ) 
 \rangle )  
 \nu^*_u ( dx ) du 
 \right] 
\\ 
\label{eq:ze}
 & = 
 \ee [ \phi ( M(s) + M^*(s) ) ]
 + 
 \frac12 
 \ee \left[ 
 \int_s^t \langle \nabla^2 \phi (M(u)+M^*(u) ) , H(u)-H^*(u) \rangle du 
 \right] 
\\
&
\nonumber
 ~~~+ 
 \ee \left[
 \int_s^t
 \int_{\real^d}
 \Psi (x, M(u) + M^*(u) )
 (\nu_u ( dx )-\nu^*_u ( dx )) du 
 \right] 
, 
\end{align} 
 where 
\begin{eqnarray*} 
 \Psi (x,y) 
 & = & 
 \phi(x+y)-\phi(y)-\sum_{i=1}^d x_i \frac{\partial \phi}{\partial y_i} ( y)
, 
 \qquad x,y\in \real^d. 
\end{eqnarray*} 
 Due to the convexity of $\phi$, 
 the Hessian $\nabla^2 \phi$ is positive semidefinite 
 hence Lemma~\ref{l31} yields 
\begin{eqnarray} 
\label{eq:ezz}
\lefteqn{ 
\! \! \! \! \! \! \! \! \! \! \! \! \! \! \! \! \! \! \! \! \! \! \! 
\! \! \! \! \! \! \! 
 \ee [ \phi ( M(t) + M^*(t) ) ] 
\leq \ee [ \phi ( M(s) + M^*(s) ) ]
} 
\\ 
\nonumber 
&& +  
 \ee \left[
 \int_s^t
 \int_{\real^d}
 \Psi (x, M(u) + M^*(u) )
 (\nu_u ( dx )-\nu^*_u ( dx )) du 
 \right] 
, 
\end{eqnarray}
 since $H^*(u)-H(u)$ is positive semidefinite 
 for fixed $(\omega,u)\in\Omega\times\real_+$. 
\\ 

\noindent 
 Finally we examine the consequences 
 of hypotheses $(i)$ and $(ii)$ on \eqref{eq:ezz}. 
\begin{description} 
\item{i)} 
 By convexity of $\phi$, $x\mapsto \Psi ( x ,y)$ is non-negative and convex 
 on $\real^d$ for all fixed $y\in \real^d$, 
 hence the second term in \eqref{eq:ezz} is non-positive. 
\item{ii)} 
 When $\nu_u$ and $\nu_u^*$ are supported by $\real_+^d$, 
 \eqref{eq:ezz} is also non-positive since 
 for all $y$, $x\mapsto \Psi (x,y)$ is non-decreasing in $x\in\real_+^d$. 
\end{description} 
 The extension to convex non ${\cal C}^2$ functions $\phi$ 
 follows by approximation of $\phi$ by an increasing sequence of ${\cal C}^2$ convex 
 Lipschitz functions, and by application of the monotone convergence theorem. 
\end{Proof} 
\begin{remark}
When $\phi\in C^2$, 
 the hypothesis on the diffusion part and on the jump part can be mixed together. 
Indeed, in order for the conclusion of Theorem \ref{11.-1} to hold 
 it suffices that 
\begin{align}
\nonumber
& \Tr(\nabla^2\phi(y) H_t)+\int_{\real^d} \Tr(\nabla^2\phi(y+\tau x)xx^t) \nu_t(dx)\\
\label{eq:iv}
&\hskip 3cm \leq\Tr(\nabla^2\phi(y) H_t^*)+\int_{\real^d} \Tr(\nabla^2\phi(y+\tau x)xx^t) \nu_t^*(dx), 
\end{align} 
 $y\in \real^d$, ${\bf 1}_{[0,1]}(\tau )d\tau dt$-a.e.
\end{remark} 
\begin{Proof} 
 Using the following version of Taylor's formula
$$
\phi (y+x) = \phi(y)+\sum_{i=1}^d x_i\phi_i'(y) +  \int_0^1 (1-\tau )\sum_{i,j=1}^d  x_ix_j\phi_{i,j}'' (y+\tau x ) d\tau, 
\qquad x,y\in \real^d
, 
$$ 
we have 
$$
 \Psi (x,y)=\int_0^1 (1-\tau )\langle \nabla^2\phi(y+\tau x ) x,x\rangle d\tau
 = \int_0^1 (1-\tau) \Tr(\nabla^2\phi(y+\tau x)xx^\dag) d\tau
$$
 and \eqref{eq:ze} rewrites as 
\begin{eqnarray*}
\lefteqn{ 
\ee[\phi(M_t+M_t^*)]-\ee[\phi(M_s+M_s^*)] 
} 
\\
 & = & 
 \frac12 
 \ee 
 \left[ \int_s^t 
 ( 
 \Tr(\nabla^2\phi(M_u+M^*_u)H_u)-\Tr(\nabla^2\phi(M_u+M^*_u)H_u^*) 
 ) 
 du\right] 
 \\
 & & + 
 \ee \left[
\int_0^1 \int_s^t
 \int_{\real^d} (1-\tau)
 \Tr(\nabla^2 \phi(M_u + M^*_u + \tau x ) xx^t)(\nu_u ( dx )-\nu^*_u ( dx )) du d\tau
 \right]
\end{eqnarray*} 
 which is non-positive from \eqref{eq:iv}.
\end{Proof}
\noindent 
 Let now $(\FF^{M}_t)_{t\in \real_+}$ and $(\FF^{M^*}_t)_{t\in \real_+}$, 
 denote the forward and backward filtrations generated by 
 $(M(t))_{t\in \real_+}$ and by $(M^*(t))_{t\in \real_+}$. 
 The proof of the following corollary of Theorem~\ref{11.-1} 
 is identical to that of Corollary~3.7 in \cite{kmp}. 
\bcor 
\label{cor1} 
 If \eqref{holds} holds and if in addition 
 $\ee [M^*(t) \mid \FF^M_t ] = 0$, $t\in \real_+$, then 
\begin{equation*} 
\label{e2} 
 \ee [\phi (M(s)+M^*(s) )] 
 \geq 
 \ee [\phi (M(t) )] 
, \qquad 0\leq s \leq t 
. 
\end{equation*} 
\ncor 
\noindent 
 In particular, if $M_0=\ee[M(t)]$ is deterministic 
 (or if $\FF_0^M$ is the trivial $\sigma$-field), 
 Corollary~\ref{cor1} shows that 
 $M(t)-\ee[M(t)]$ is more convex concentrated than $M^*_0$, 
 i.e.: 
$$\ee [\phi (M(t)-\ee[M(t)])] \leq \ee [\phi (M^*_0)], 
 \qquad t\geq 0, 
$$ 
 for all sufficiently integrable convex functions $\phi$ on 
 $\real^d$. 
\noindent 
 In applications to convex concentration inequalities 
 the independence of $(M(t))_{t\in \real_+}$ with $(M^*(t))_{t\in \real_+}$ 
 will not be required, see Section~\ref{sc01}. 
\\ 

\noindent 
 Note that in case $\nu^* (dt,dx )$ has the form 
\begin{equation*} 
\label{car3} 
\nu^* (dt,dx ) = \lambda^*(t) \delta_{k}(dx) dt 
, 
\end{equation*} 
 where $k\in \real^d$ and $(\lambda^*(t))_{t\in \real_+}$ is a positive $\FF^*_t$-predictable process, 
 then condition $(i)$ (resp. $(ii)$) of 
 Theorem~\ref{11.-1} is equivalent to: 
$$
\nu_t = \lambda (t) \delta_k  
\qquad  \mbox{and}  \qquad 
 \lambda(t)\leq \lambda^*(t) 
$$ 
resp. to: $k\in(\real_+)^d$, $\nu_t(\real^d) \leq \lambda^*(t)$ and
$$
\nu_t 
 \left( 
 \real^d \setminus 
 \bigcap_{i=1}^d ]-\infty, k_i] 
 \right) 
 = 
 0,
$$
 i.e. the jump $\Delta M_i(t)$ is a.s. upper bounded by $k_i$, 
 $i=1,\ldots ,d$. 
\\

\noindent 
 Theorem \ref{11.-1} applies for instance when the jump parts of
 $(M(t))_{t\in \real_+}$ and of $(M^*(t))_{t\in \real_+}$ are point processes. 
 Let $(W(t))_{t\in \real_+}$ be a standard $\real^n$-valued Brownian motion and 
 $(W^*(t))_{t\in \real_+}$ be a backward standard $\real^n$-valued Brownian motion, 
 and let $(Z(t))_{t\in \real_+}$ be a point process in $\real^n$ given by 
 $Z(t)=(Z_1 (t), \dots, Z_n(t))$ where $(Z_i(t))_{t\in \real_+}$ 
 is a real point process with intensity $(\lambda_i (t))_{t\in \real_+}$, 
 $1\leq i\leq n$. 
 Similarly, let $(Z^*(t))_{t\in \real_+}$ be a backward point process 
 in $\real^n$ with intensity 
 $\lambda^*(t)=(\lambda^*_1(t), \dots, \lambda^*_n(t))$, $t\in\real_+$. 
\\ 

\noindent 
 We can take 
\begin{equation*} 
\label{dec1.0} 
 M(t) = M_0 + \int_0^t A(s) dW(s) + \int_0^t J(s) (dZ(s) - \lambda(s) ds), \quad t\in \real_+, 
\end{equation*} 
 and 
\begin{equation*} 
\label{dec2.0} 
 M^*(t) = \int_t^{+\infty} A^*(s) d^*W^*(s) + \int_t^{+\infty} J^*(s) (d^* Z^*(s) - \lambda^*(s) ds), 
 \quad t\in \real_+, 
\end{equation*} 
 where $(A(t))_{t\in \real_+ }$, $(J(t))_{t \in \real_+ }$, 
 resp. $(A^*(t))_{t\in \real_+ }$, $(J^*(t))_{t \in \real_+ }$ 
 are ${\cal M}^{d\times n}$-valued and predictable with respect to 
$$\FF^M_t:=\sigma( W(s), Z(s) : s\leq t) 
, 
 \qquad 
 \mbox{resp.} 
 \qquad 
 \FF^{M^*}_t:= \sigma( W^*(s), Z^*(s) : s\geq t), 
$$ 
 $t\in \real_+$, i.e. 
$$ 
 M_i(t)= M_i(0) + \sum_{j=1}^n 
 \int_0^t A_{i,j}(s) dW_j(s) + \sum_{j=1}^n \int_0^t J_{i,j}(s) (dZ_j(s) - \lambda_j(s) ds) 
$$
 and 
$$
 M^*_i(t)= M^*_i(0) + \sum_{j=1}^n \int_t^\infty 
 A_{i,j}(s) dW^*_j(s) + \sum_{j=1}^n \int_t^\infty 
 J_{i,j}(s) (d^*Z^*_j(s) - \lambda_j(s) ds) 
, 
$$ 
 $t\in \real_+$, $i=1,\ldots , d$, with 
\begin{equation} 
\label{car2a}
 \nu_t ( dx ) = 
 \sum_{j=1}^n \lambda_j(t) \delta_{(J_{1,j}(t),\dots,J_{d,j}(t))} (dx) 
\end{equation} 
 and 
\begin{equation} 
\label{car2b}
 \nu^*_t ( dx ) = 
\sum_{j=1}^n \lambda^*_j(t)  \delta_{(J^*_{1,j}(t), \dots,J^*_{d,j}(t))} (dx). 
\end{equation} 
 As seen above, condition $(i)$ and $(ii)$ of Theorem~\ref{11.-1} imply 
$$
 \sum_{j=1}^n \lambda_j(t) J_{k,j}(t) 
 \leq 
 \sum_{j=1}^n \lambda^*_j(t) J^*_{k,j}(t) 
, 
 \qquad 
 k=1,\ldots , d, 
$$ 
 and under both conditions we have 
$$
 \sum_{j=1}^n \lambda_j(t) 
 \leq 
 \sum_{j=1}^n \lambda^*_j(t) 
. 
$$ 
 More details will be given in Section~\ref{s5} 
 on the meaning of conditions $(i)$ and $(ii)$ 
 of Theorem~\ref{11.-1} imposed on $\nu_t$ and $\nu^*_t$ 
 defined in \eqref{car2a} and \eqref{car2b}
 for the order $\preceq_{\rm cxp}$. 
\\ 

\noindent 
 Conditions \eqref{hp1} \eqref{hp2} will hold in particular when 
$$
\FF_t=\FF_t^M\vee\FF_0^{M^*}\; \mbox{ and } \;
\FF_t^*=\FF_\infty^{M}\vee\FF_t^{M^*}, 
 \qquad 
 t\in \real_+
, 
$$
 see Section~\ref{sc01}. 
\section{Convex ordering and predictable representation} 
\label{sc01} 
 Let $(W(t))_{t\in \real_+}$ be a $n$-dimensional Brownian motions and 
 $\mu (dx,dt)$ be a jump measure with jump characteristics of the form 
\begin{equation} 
\label{aa1}
\nu(dt,dx)=\nu_t(dx) dt 
, 
\end{equation} 
 generating a filtration $(\FF^M_t)_{t\in \real_+}$. 
\\ 

\noindent 
 Consider $(A(t))_{t\in \real_+}$ and $( \hat{A} (t))_{t\in \real_+}$ 
 two ${\cal M}^{d\times n}$-valued, $\FF^M_t$-predictable 
 square-integrable processes, and $(t,x)\mapsto B_t(x)$ and 
 $(t,x)\mapsto \hat{B}_t(x)$ two $\real^d$-valued
 $\FF^M_t$-predictable processes in 
 $L^1(\Omega\times \real^d\times\real_+, dP\nu_t(dx)dt)$. 
\bthm 
\label{112} 
 Let $(\tilde{W}(t))_{t\in \real_+}$ be an $n$-dimensional Brownian motion and 
 $\tilde{\mu} (dx,dt)$ be a jump measure with jump characteristic of the form 
 $\tilde\nu(dt,dx)=\tilde\nu_t(dx) dt$, 
 both independent of ${\cal F}^M$, and consider 
$$ 
  F = \int_0^\infty A( t ) dW(t) +\int_0^\infty  \int_{\real^n} B_t (x) (\mu(dt ,dx) -\nu_t (dx)dt)
$$ 
 and 
$$ 
 G = \int_0^\infty \hat{A} (t) d\tilde W(t) +\int_0^\infty \int_{\real^n } \hat{B}_t (x) (\tilde\mu(dt ,dx) -\tilde\nu_t (dx)dt )
. 
$$ 
 Assume that 
$$ 
 A^\dag (t)A(t)\leq_{\rm psd} \hat{A}^\dag (t)\hat{A} (t), 
\qquad 
 dPdt-a.e. 
, 
$$ 
 and that for almost all $t\in \real_+$, we have either: 
\begin{description} 
\item{i)} 
 $\nu_t \circ B_t^{-1} \preceq_{\rm cxp} \tilde\nu_t \circ \hat{B}_t^{-1}$, 
 $P$-a.s., 
\item or: 
\item{ii)} $B_t$ and $\hat{B}_t$ are non-negative and 
 $\nu_t \circ B_t^{-1} \preceq_{\rm cxpi} \tilde\nu_t \circ \hat{B}_t^{-1}$, 
 $P$-a.s. 
\end{description} 
\noindent Then for all 
 convex functions $\phi : \real^d \to \real$ we have 
\begin{equation}
\label{hold4}
 \ee[\phi( F )]\leq \ee[\phi( G )]
. 
\end{equation}
\nthm 
\bprf 
 Again we start by assuming that $\phi$ is a Lipschitz convex function. 
 Let $(M(t))_{t \in \real_+}$ denote the forward martingale 
 defined as 
$$ 
 M(t) = \int_0^t A(s) dW(s) +\int_0^t  \int_{\real^n} B_s(x) (\mu(ds,dx) -\nu_s(dx)ds), 
$$ 
 $t\in \real_+$, 
 let $(\FF^{\tilde M}_t)_{t\in \real_+}$ 
 denote the backward filtration generated by 
 $\{ \tilde{W} (t) , \tilde{\mu} (dx,dt) \}$, 
 and let 
$$ 
\FF_t=\FF_t^M\vee \FF_\infty^{\tilde M} 
 \qquad 
 \mbox{and} 
 \qquad 
 \FF_t^* = 
 \FF_\infty^M \vee 
 \FF^{\tilde M}_t 
, 
$$ 
 so that $(M(t))_{t\in \real_+}$ is an $\FF_t$-forward martingale. 
 Since $(\hat{A} (t) )_{t\in \real_+}$, $(\hat{B}_t)_{t\in \real_+}$ 
 are $\FF^M_t$-predictable, 
 the processes $(\hat{A} (t) )_{t\in \real_+ }$ and 
 $(\hat{B}_t)_{t\in \real_+ }$ are 
 independent of $(\tilde W_t)_{t\in \real_+}$ and of $\tilde \mu (dt,dx)$. 
 In this case, the forward and backward differentials 
 coincide and the process $(M^*(t))_{t\in \real_+ }$ defined as 
$$ 
M^*(t) =\ee[ G \mid \FF_t^*]=\int_t^\infty  \hat{A} (s) d\tilde W(s) +\int_t^\infty 
 \int_{\real^n} \hat{B}_s(x) (\tilde\mu(ds,dx) -\tilde\nu_s(dx)ds)
, 
$$ 
 $t \in \real_+$, is an $\FF_t^*$- backward martingale with 
 $M^*(0)=G$. 
%
 Moreover the jump characteristics of $(M_t)_{t\in \real_+}$ 
 and of $(M^*_t)_{t\in \real_+}$ are 
$$ 
 \nu_M (dx) = {\bf 1}_{\real^d\setminus \{ 0 \} } (x) 
 \nu_t \circ B_t^{-1}(dx) 
 \qquad \mbox{ and } \qquad 
 \nu_{M^*} (dx) = 
 {\bf 1}_{\real^d\setminus \{ 0 \} } (x) 
 \tilde\nu_t \circ \hat{B}_t^{-1}(dx) 
.
$$
 Applying Theorem~\ref{11.-1} to the forward and 
 backward martingales $(M(t))_{t\in \real_+}$ 
 and $(M^*(t))_{t\in \real_+}$ yields 
$$
 \ee[\phi(M(t)+M^*(t))]\leq\ee[\phi(M(s)+M^*(s))], 
 \qquad 
 0\leq s\leq t, 
$$ 
 for all convex functions $\phi$ and for $0\leq s\leq t$.
 Since $M(0)=0$, $M^*(0)=G$ and $\lim_{t\to \infty} M^*(t)= 0$ in $L^2(\Omega )$, 
 we obtain \eqref{hold4} for convex Lipschitz function $\phi$  by 
 taking $s=0$ and letting $t$ go to infinity. 
 Finally we extend the formula to all convex integrable 
 functions $\phi$ by considering 
 an increasing sequence of Lipschitz convex functions $\phi_n$ 
 converging pointwise to $\phi$. 
 Applying the monotone convergence theorem to the non-negative sequence 
 $\phi_n(F)-\phi_0(F)$, we have
 $$
 \ee[\phi(F)-\phi_0(F)] =\lim_{n\to \infty}\ee[\phi_n(F)-\phi_0(F)]
, 
 $$
 which yields $\ee[\phi(F)] =\lim_{n\to \infty}\ee[\phi_n(F)]$ 
 since $\phi_0(F)$ is integrable. 
 We proceed similarly for $\phi(G)$, allowing us 
 to extend \eqref{hold4} to the general case. 
\nprf 
\noindent 
 Note that if $(\hat{A} (t) )_{t\in \real_+ }$ and 
 $(\hat{B}_t)_{t\in \real_+ }$ are deterministic 
 then $(\tilde{W}(t))_{t\in \real_+}$ and $\tilde{\mu} (dx,dt)$ 
 can be taken equal to $(W(t))_{t\in \real_+}$ and $\mu (dx,dt)$ respectively. 
\subsubsection*{Example: point processes}
 Let $(A(t))_{t\in \real_+}$, 
 $( \hat{A} (t))_{t\in \real_+}$, 
 $( W (t))_{t\in \real_+}$ and 
 $(\tilde{W}(t))_{t\in \real_+}$ 
 be as in Theorem~\ref{112} above and consider 
$$ 
 Z(t)=(Z_1(t), \dots, Z_n(t)) 
 \quad 
 \mbox{and} 
 \quad 
 \tilde Z(t)=(\tilde Z_1(t), \dots, \tilde Z_n(t)) 
$$ 
 to be two independent point processes in $\real^n$ with 
 compensators 
$$
\sum_{i=1}^n \lambda_i(t)\delta_{e_i}
\quad \mbox{ and } \quad
\sum_{i=1}^n \tilde\lambda_i(t)\delta_{e_i} 
, 
$$
 where $e_i = (0,\dots,1,\dots,0)$, $i=1,\ldots ,n$, 
 denotes the canonical basis in $\real^n$, and let 
$$
 \FF_t^M = \sigma ( 
 W(s) , \ Z(s) \ : \ 0 \leq s\leq t ) 
, 
 \qquad 
 t\in \real_+ 
. 
$$ 
\bcor 
 Given $J(t)$ and $\hat{J} (t)$ two ${\cal M}^{d\times n}$-valued integrable 
 ${\cal F}^M_t$-predictable processes, let 
$$ 
 F =\int_0^\infty A( t ) dW( t ) + \int_0^\infty J(t) (dZ(t)-\lambda(t)dt) 
, 
$$ 
 and 
$$ 
 G =\int_0^\infty \hat{A} ( t ) d\tilde W( t )+\int_0^\infty 
 \hat{J} ( t ) (d\tilde Z( t )-\tilde\lambda( t )d t ). 
$$ 
 Assume that 
$$ 
 A^\dag (t)A(t)\leq_{\rm psd} \hat{A}^\dag (t)\hat{A} (t), 
\qquad 
 dPdt-a.e. 
, 
$$ 
 and that for almost all $t\in \real_+$, we have either: 
\begin{description}
\item{i)} 
 $ 
 \sum_{j=1}^n\lambda_j(t) \delta_{(J_{1,j}(t),\dots,J_{d,j}(t))} 
 \preceq_{\rm cxp} 
 \sum_{j=1}^n \tilde\lambda_j(t)  \delta_{(\hat{J}_{1,j}(t), \dots,\hat{J}_{d,j}(t))} 
 $, $P$-a.s.,  
\item or: 
\item{ii)} 
 $J_{i,j}(t)\geq 0$, and $\hat{J}_{i,j}(t)\geq 0$, $i=1,\ldots ,d$, $j=1,\ldots ,n$, 
 and 
$$ 
 \sum_{j=1}^n\lambda_j(t) \delta_{(J_{1,j}(t),\dots,J_{d,j}(t))} 
 \preceq_{\rm cxpi} 
 \sum_{j=1}^n \tilde\lambda_j(t)  \delta_{(\hat{J}_{1,j}(t), \dots,\hat{J}_{d,j}(t))} 
, \quad P\mbox{-a.s}. 
$$ 
\end{description} 
\noindent Then for all 
 convex functions $\phi : \real^d \to \real$ we have 
$$ 
 \ee[\phi( F )]\leq \ee[\phi( G )]
. 
$$
\ncor 
\bprf 
 We apply Theorem \ref{112} with 
 $B_t(x):=J(t) x$ and $\hat{B}_t(x):=\hat{J}(t) x$, $t\in \real_+$, 
 and 
\begin{equation} 
\label{car3.1} 
 \nu_t (dx) 
 =
 {\bf 1}_{\real^d\setminus \{ 0 \} } (x) 
 \sum_{j=1}^n\lambda_j(t) \delta_{(J_{1,j}(t),\dots,J_{d,j}(t))}  
 \circ B_t^{-1} (dx) 
\end{equation} 
 and 
\begin{equation} 
\label{car3.1.2} 
 \tilde \nu_t (dx) 
 = 
 {\bf 1}_{\real^d\setminus \{ 0 \} } (x) 
 \sum_{j=1}^n \tilde\lambda_j(t)  \delta_{(\hat{J}_{1,j}(t), \dots,\hat{J}_{d,j}(t))} 
 \circ \hat{B}_t^{-1} (dx) 
.
\end{equation} 
\nprf
\noindent 
 Note that $\nu_t \preceq_{\rm cxp} \tilde \nu_t$ if and only if 
 $\lambda_j(t) \leq \tilde\lambda_j(t)$, $j=1,\ldots ,n$, since 
 $\Supp (\nu_t)=\Supp (\mu_t) = \{ e_1,\ldots ,e_n\}$. 
 In Section~\ref{s5} we will give a geometric interpretation 
 of the convex ordering $\preceq_{\rm cxp}$, with application 
 to the conditions imposed on $\nu_t \circ B_t^{-1}$ 
 and $\tilde \nu_t \circ \hat{B}_t^{-1}$ in $(i)$ and $(ii)$ above. 
\subsubsection*{Example: Poisson random measures} 
 Consider $\sigma$, $\tilde{\sigma}$ two atomless Radon measures 
 on $\real^n$ with 
$$ 
\int_{\real^n} ( | x  |^2 \wedge 1 ) \sigma (dx) < \infty,
\quad 
 \mbox{and} 
\quad 
\int_{\real^n} ( | x  |^2 \wedge 1 ) \tilde{\sigma} (dx) < \infty 
, 
$$ 
 and two Poisson random measures 
$$
 \omega (dt , dx ) = \sum_{i\in \nn} \delta_{(t_i,x_i)}(dt , dx) 
 \quad 
 \mbox{and} 
 \quad 
\tilde{\omega} (dt , dx ) = \sum_{i\in \nn} \delta_{(\tilde{t}_i,\tilde{x}_i)}(dt , dx)
$$ 
 with respective intensities $\sigma (dx) dt$ and $\tilde{\sigma} (dx) dt$ 
 on $\real^n \times \real_+$ under $P$. 
 Let also $(W(t))_{t\in \real_+}$ and $(\tilde W(t))_{t\in\real_+}$ 
 be independent $n$-dimensio\-nal standard Brownian motions, 
 independent of 
 $\tilde{\omega} (dt,dx)$ 
 under $P$ and let $(A(t))_{t\in \real_+}$, $( \hat{A} (t))_{t\in \real_+}$ 
 be as in Theorem~\ref{112} above, with 
$$
\FF_t^M = \sigma ( 
 W(s) , \ \omega ([0,s]\times A) \ : \ 0 \leq s\leq t, \ A\in {\cal B}_b (\real^n )) 
, 
 \qquad 
 t\in \real_+, 
$$ 
 where ${\cal B}_b (\real^n ) = \{ A\in {\cal B} (\real^n ) \ : \ \sigma (A) < \infty \}$. 
\bcor 
\label{bcor} 
 Let $(J_{t,x})_{(t,x) \in \real_+\times \real^n}$ be a 
 $\real^d$-valued $\FF_t^M$-predictable 
 process, integrable with respect to $dPdt\sigma(dx)$, 
 and let $(\hat{J}_{t,x})_{(t,x)\in \real_+\times \real^n}$ be 
 an $\real^d$-valued deterministic function, 
 integrable with respect to $dt\tilde \sigma(dx)$. 
 Consider the random variables 
\begin{equation} 
\label{erf} 
 F = \int_0^\infty A(t) dW(t) 
 +\int_0^\infty\int_{\real^n} J_{t^-,x}(\omega(dt,dx)-\sigma(dx)dt)
\end{equation} 
 and 
$$ 
 G = \int_0^\infty \hat{A} (t) d\tilde W(t) 
 +\int_0^\infty \int_{\real^n} \hat{J}_{t^-,x}(\tilde \omega(dt,dx)-\tilde\sigma(dx)dt)
. 
$$ 
 Assume that 
$$ 
 A^\dag (t)A(t)\leq_{\rm psd} \hat{A}^\dag (t)\hat{A}(t), 
\qquad 
 dPdt\mbox{-a.e.} 
, 
$$ 
 and that for almost all $t\in \real_+$, we have either: 
\begin{description}
\item{i)} 
 $ 
 \sigma \circ J_{t^-,\cdot}^{-1} 
 \preceq_{\rm cxp} 
 \tilde\sigma \circ \hat{J}_{t^-,\cdot}^{-1}
 $, $P$\mbox{-a.s.}, 
\item or: 
\item{ii)} 
 $\hat{J}_{t^-,x} \geq 0$, $\sigma (dx)$\mbox{-a.e.}, and 
$$ 
 \sigma \circ J_{t^-,\cdot}^{-1} 
 \preceq_{\rm cxpi} 
 \tilde\sigma \circ \hat{J}_{t^-,\cdot}^{-1}, 
 \quad 
 P\mbox{-a.s}. 
$$ 
\end{description} 
\noindent Then for all 
 convex functions $\phi : \real^d \to \real$ we have 
$$ 
 \ee[\phi( F )]\leq \ee[\phi( G )]
. 
$$
\ncor 
\bprf 
 We apply Theorem~\ref{112} with the jump characteristics 
$$\nu_t (dx) =  {\bf 1}_{\real^d\setminus \{ 0 \} } (x) 
 \sigma \circ J_{t^-,\cdot}^{-1} (dx) 
 \qquad 
 \mbox{and} 
 \qquad 
 \tilde\nu_t (dx) = 
 {\bf 1}_{\real^d\setminus \{ 0 \} } (x) 
 \tilde\sigma \circ \hat{J}_{t^-,\cdot}^{-1} (dx) 
, 
$$ 
 and $B_t(x)=\hat{B}_t(x)=x$, $x\in \real^n$. 
\nprf 
\noindent 
 Condition $(i)$, resp. $(ii)$ in Corollary~\ref{bcor} can be written as 
$$
\int_{\real^n} f(J_{t^-,x}) \sigma(dx)\leq \int_{\real^n} f(\hat{J}_{t^-,x}) \tilde\sigma(dx)
$$
 for all non-negative convex functions $f:\real^d\to\real$, 
 resp. for all non-negative and non-decreasing convex 
 functions $f:\real^d\to\real$.
 In particular, Corollary~\ref{bcor}-$ii)$ holds if  we have 
 $\sigma\preceq_{\rm cxpi}\tilde\sigma$ and 
 $J_{t,x}\leq \hat{J}_{t,x}$, $dt\sigma (dx)dP$-a.e., and if 
 $x \mapsto J_{t,x}$, $x\mapsto \hat{J}_{t,x}$ 
 are non-decreasing and convex on $\real^n$ for all $t\in \real_+$. 
\\ 

\noindent 
 We may also apply Theorem \ref{112} to $F$ as in \eqref{erf} with 
$$
 G =\int_0^\infty \hat{A} (t) d\tilde W(t)+\int_0^\infty \hat{J} (t) 
 (d\tilde Z(t)-\tilde\lambda(t)dt)
$$
 where $(\hat{A} (t))_{t\in\real_+}$ and $(\hat{J} (t) )_{t\in\real_+}$ are 
 $\FF^M_t$-predictable ${\cal M}^{d\times n}$-valued processes 
 and $(\tilde Z_1, \dots, \tilde Z_n)$ is a $\real^n$-valued point process 
 independent of ${\cal F}^M$, with $\nu_t=\sigma \circ J_{t^-,\cdot}^{-1}$ 
 and 
$$ 
 \tilde\nu_t=\sum_{i=1}^n \tilde\lambda_i(t)\delta_{e_i} 
. 
$$ 
 In case $x\mapsto J_{t,x}$ is convex (resp. non-negative, non-decreasing and convex) 
 on $\real^n$ for all $0\leq t\leq T$, condition $(i)$ resp. $(ii)$, 
 of Theorem~\ref{112} is satisfied provided 
\begin{equation} 
\label{car5} 
 \sigma \preceq_{\rm cxp}\sum_{i=1}^n \tilde\lambda_i(t)\delta_{e_i}, 
 \qquad 
 \mbox{resp.} 
 \qquad 
 \sigma \preceq_{\rm cxpi}\sum_{i=1}^n \tilde\lambda_i(t)\delta_{e_i}. 
\end{equation} 
\section{A geometric interpretation for discrete measures} 
\label{s5} 
 The next lemma provides a first interpretation of the 
 order $\preceq_{\rm cxp}$. 
\blem
 If $\mu$ and $\nu$ are two measures on $\real^d$ with finite 
 supports, then $\mu \preceq_{\rm cxp} \nu$ implies 
\begin{equation} 
\label{cnd00} 
{\mathscr C} (\Supp (\mu ) ) \subset {\mathscr C} ( \Supp (\nu ) ), 
\end{equation} 
 where ${\mathscr C} (A)$ denote the convex hull of any subset 
 $A$ of $\real^d$. 
\nlem
\begin{Proof} 
 Let $H$ be any half-space of $\real^d$ such that 
 $\Supp(\nu) \subset H$. 
 For any convex function $\phi$ such that $\{\phi\leq 0\}=H$ 
 and $\phi_{|\partial H}=0$ we have 
 $\int_{\real^d} \phi^+ \nu(dx)=0$, hence 
 $\int_{\real^d} \phi^+ \mu(dx)=0$ 
 since $\mu \preceq_{\rm cxp} \nu$, 
 which implies $\Supp(\mu) \subset H$. 
 The conclusion follows from the characterization of 
 the convex hull ${\mathscr C} (\Supp (\mu ) )$, 
 resp. ${\mathscr C} ( \Supp (\nu ) )$, 
 as the intersections of all half-spaces containing it. 
\end{Proof} 

\noindent 
 However the necessary condition \eqref{cnd00} 
 is clearly not sufficient to ensure the convex ordering of $\mu$ and $\nu$. 
 Our aim in this section is to find a more precise geometric interpretation of 
 $\mu \preceq_{\rm cxp} \nu$ in the 
 case of finite supports, with the aim of applying this criterion 
 to the jump measures defined in \eqref{car2a}, \eqref{car2b}, \eqref{car3.1}, \eqref{car3.1.2}  
 and \eqref{car5}. 
\\ 

\noindent 
 For all  $u\in S^{d-1}$ the unit sphere in $\real^d$, let 
 $\mu_u = \mu\circ \langle u,\cdot\rangle^{-1}$ 
 (resp. $\nu_u = \nu\circ \langle u,\cdot\rangle^{-1}$) 
 denote the image of $\mu$, resp. $\nu$, on $\real$ by the mapping 
 $x\mapsto \langle u, x\rangle$. 
 We have $\mu_u \preceq_{\rm cxp} \nu_u$ and 
 the survival function $\phi_{\mu,u}$  associated with $\mu_u$, 
 defined by 
$$
 \phi_{\mu,u}(a) 
 = 
 \int_\real (y-a)^+\,d\mu_u(y) 
 = 
 \int_{\real^d} (\langle y,u\rangle -a)^+\,d\mu(y) 
, 
$$
 is a convex function with $\phi_{\mu,u}\le \phi_{\nu,u}$ for all $u\in S^{d-1}$. 
 Moreover for all $a\in \real$ such that $a$ is sufficiently large we have 
 $\phi_{\mu,u}(a)=\phi_{\nu,u}(a) = 0$.
\\ 

\noindent 
For every $x\in \real^d$ and $u\in S^{d-1}$ let 
\begin{equation}
\label{E1}
a_{x,u} =\inf\left\{ b\in \real \ : \ b \ge \langle u,x\rangle,\ \phi_{\mu,u}(b)=\phi_{\nu,u}(b)\right\}
\end{equation}
and 
$$
{\cal D}_{x,u}= \{y\in \real^d \ : \ \langle u,y\rangle \le a_{x,u} \} 
, 
$$ 
 which is a closed half-space containing $x$. Finally we let 
$$
{\cal C}_x : =\bigcap_{u\in S^{d-1}}{\cal D}_{x,u} 
$$ 
 which is a compact convex set containing $x$. 
 On the other  hand, letting 
$$ 
 \tilde{{\cal D}}_u = \{ z\in \real^d \ : \ \langle u,z\rangle \le \tilde{a}_u \} 
, 
\qquad x\in \real^d, \quad u\in S^{d-1}, 
$$ 
where
$$ 
 \tilde{a}_u = 
 \inf\left\{ b\in \real \ : \ b \ge \langle u,y\rangle,\ \forall y\in \Supp (\nu ) 
 \right\}
, 
$$ 
 we have 
$$ 
 {\mathscr C} 
 ( \Supp (\nu ) ) 
 = 
 \bigcap_{\it u\in S^{d-1}} 
 \tilde{{\cal D}}_{\it u} 
. 
$$ 
 Note that we have ${\cal C}_x \subset {\mathscr C} ( \Supp (\nu ) )$ 
 if $\mu \preceq_{\rm cxp} \nu$, indeed 
 we have  ${\cal D}_{x,u} \subset \tilde{{\cal D}}_u$ 
 since $a_{x,u}\leq \tilde{a}_u$, as follows from 
 $\phi_{\mu ,u}(b) = \phi_{\nu ,u} (b) =0$ for all $b\geq \tilde{a}_u$. 
 On the other hand, if $\mu = \delta_x \preceq_{\rm cxp} \nu$ 
 then 
\begin{equation} 
\label{c}  
 {\cal C}_x = {\mathscr C} ( \Supp (\nu)) 
\end{equation} 
 since for all $u\in S^{d-1}$ 
 there exists $z\in \Supp (\nu )$ such that 
 $\langle u,z\rangle = \tilde{a}_u$ and for all 
 $b\in  ( \langle u,x\rangle , \tilde{a}_u]$ we have 
$$ 
 \phi_{\nu , u} (b) \geq (\tilde{a}_u - b ) 
 \nu ( \{ z \})>0 
$$ 
 and $\phi_{\mu ,u}(b)=0$, implying $a_{x,u}=\tilde{a}_u$ 
 and ${\cal D}_{x,u} = \tilde{{\cal D}}_u$. 
%
 Note that in Theorem~\ref{T5.1} below 
 the existence of $x\in \real^d$ such that $\mu(\{x\})>\nu(\{x\})$ 
 is always satisfied when $\mu \preceq_{\rm cx} \nu$ and $\mu \not= \nu$. 
\bthm
\label{T5.1} 
 Assume that $\mu$ and $\nu$ have finite supports and that 
 $\mu \preceq_{\rm cxp} \nu$. 
 Then for all 
 $x\in \real^d$ such that $\mu(\{x\})>\nu(\{x\})$ 
 there exists 
 $k\in \{ 2, \ldots , d+1\}$ 
 and $k$ elements $\displaystyle y_1,\ldots,y_k \in \Supp (\nu)$ 
 distinct from $x$, such that 
\begin{equation} 
\label{cnd01} 
 \{ x \} 
 \ 
 \subset 
 \ 
 {\mathscr C} (\{y_1,\ldots,y_k\}) 
 \ 
 \subset 
 \ 
 {\cal C}_x 
. 
\end{equation}
\nthm
\bprf
 We only need to prove that $x$ belongs to the convex hull of
$({\cal C}_x\backslash\{x\})\cap \,\Supp (\nu)$. 
 Indeed, if $x\in {\mathscr C} ( ({\cal C}_x\backslash\{x\})\cap \,\Supp (\nu))$ 
 then there exists $k$ points $y_1,\ldots ,y_k$ 
 in this set, $k\geq 2$, 
 such that $x$ is the convex barycenter of $y_1,\ldots ,y_k$, and the 
 Caratheodory theorem (see e.g. \cite{rockafeller}, Theorem~17.1) 
 shows that the conclusion holds 
 for some $k\in \{2,\ldots,d+1\}$. 
\\ 

\noindent 
 Assume now that the assertion of the theorem is true when $\mu$ and $\nu$
have disjoint supports, and let $\mu$ and $\nu$ be any measures with
finite supports, such that $\mu \preceq_{\rm cxp} \nu$ and $\mu\not=\nu$. We
let ${\cal S}^+=\{y\in \real^d,\ \nu(\{y\})\ge\mu(\{y\})\}$,  ${\cal S}^-=\{y\in \real^d,\
\mu(\{y\})>\nu(\{y\})\}$. 
 These sets are not empty since $\mu\not=\nu$. 
 Let $\mu'$ and $\nu'$ the measures defined by 
$$
\mu'(\{y\})=\mu(\{y\})-\nu(\{y\}),\qquad \nu'(\{y\})=0, \qquad y\in {\cal S}^-,
$$
$$
\mu'(\{y\})=0,\qquad \nu'(\{y\})=\nu(\{y\})-\mu(\{y\}), \qquad y\in {\cal S}^+.
$$
 If $f$ is a function on $\real^d$ we have 
$$
\mu'(f)=\mu(f)-\sum_{y\in {\cal S}^-}\nu(\{y\})f(y)-\sum_{y\in {\cal S}^+}\mu(\{y\})f(y)
$$
and 
$$
\nu'(f)=\nu(f)-\sum_{y\in {\cal S}^-}\nu(\{y\})f(y)-\sum_{y\in {\cal S}^+}\mu(\{y\})f(y)
$$
which implies that $\mu \preceq_{\rm cxp} \nu$ if and only
if $\mu' \preceq_{\rm cxp} \nu'$. It also implies that for
all $u\in S^{d-1}$ and $b\in \real$, 
$$
\phi_{\mu,u}(b)=\phi_{\nu,u}(b)\Longleftrightarrow \phi_{\mu',u}(b)=\phi_{\nu',u}(b).
$$
 From this together with the fact that $\phi_{\mu',u}\le \phi_{\nu',u}$, we conclude that 
 if ${\cal D}'_{x,u}$ is defined as ${\cal D}_{x,u}$ 
 but with $(\mu,\nu)$ replaced
by $(\mu',\nu')$, then ${\cal D}'_{x,u}={\cal D}_{x,u}$. 
 Finally remarking that the support
of $\nu'$ is included in the support of $\nu$, we proved that it is
sufficient to do the proof with $\mu'$ and $\nu'$. 
\\ 

\noindent 
So in the sequel we assume that $\mu$ and $\nu$ have disjoint supports. 
As a consequence of Theorem~40 in \cite{delmey2} applied 
 to the cone of non-negative convex functions, 
 there exists an admissible\footnote{Admissible means here that for every $x\in \real^d$ we have 
 $\delta_x \preceq_{\rm cxp} K(x,dy)$.} 
 kernel $K$ such that $\mu K=\nu$. 
\\ 
 
\noindent 
 Let now $x\in \real^d$ satisfy  $\mu(\{x\})>\nu(\{x\})=0$. 
 Clearly the support of $K(x,dy)$ is included in the
support of $\nu$, and by \eqref{c}, 
 $x$ is in the convex hull of the support of $K(x,dy)$. 
 Finally we are left to prove that the
support of $K(x,dy)$ is included in ${\cal C}_x$. 
\\ 

\noindent 
 For this we let $\mu^x$ be the 
 measure defined by
$$\mu^x(\{y\})=\left\{ 
 \begin{array}{ll} 
 \mu(\{y\})+\mu(\{x\})K(x,\{y\}), & \qquad y\not=x, 
\\ 
\\ 
0, & \qquad y=x 
. 
\end{array} 
\right. 
$$
 Then  $\mu \preceq_{\rm cxp} \mu^x$ and $\mu^x \preceq_{\rm cxp} \nu$, 
 which is easily proved by the existence of
admissible  kernels $P$ and $P'$ such that $\mu P=\mu^x$ and $\mu^x
P'=\nu$. 
 More precisely they are  given by  
$$ 
P(x, dy)=\sum_{z\in 
 \Supp (\nu)}K(x,\{z\})\delta_z(dy), \qquad P(z,dy) =\delta_z(dy),\quad z\not=x, 
$$ 
 and 
$$
P'(x',dy) = 
 \left\{ 
 \begin{array}{ll} 
 \displaystyle \sum_{z\in \Supp (\nu)}K(x',\{z\})\delta_z(dy), & x'\in  
 \Supp (\mu)\backslash\{x\}, 
\\ 
\\ 
 \delta_{x'}(dy), & x'\not\in 
 \Supp (\mu)\backslash\{x\}. 
\end{array} 
\right. 
$$ 
 So for every $u\in
S^{d-1}$, we have 
$\displaystyle
\phi_{\mu,u}\le \phi_{\mu^x,u}\le\phi_{\nu, u}.
$
Let us prove that this inequality implies that any point of the
support of $K(x,dy)$ belongs to ${\cal D}_{x,u}$. Assume that a point $z$  of the
support of  $K(x,dy)$ does not. An easy calculation shows that the
right derivatives $\phi_{\mu^x,u}'$ and $\phi_{\mu,u}'$ satisfy
$\phi_{\mu^x,u}'(t)=-\mu_u^x(]t,\infty[)$ and
    $\phi_{\mu,u}'(t)=-\mu_u(]t,\infty[)$. From the definition of
    $\mu^x$ we see that for $t\not=\langle u,x\rangle$,
    $\mu_u^x(\{t\})\ge \mu_u(\{t\})$. This implies that for $t\ge
    \langle u,x\rangle$, $$-\mu_u^x(]t,\infty[)\le
    -\mu_u(]t,\infty[).$$ But since $\mu_u^x(\{\langle
    u,z\rangle\})>\mu_u(\{\langle u,z\rangle\})$ and $\langle
    u,z\rangle >\langle u,x\rangle$, we have for $t\in
    [\langle u,x\rangle,\langle u,z\rangle[$ 
$$
-\mu_u^x(]t,\infty[)<
    -\mu_u(]t,\infty[)
$$
and this implies that $\phi_{\mu,u}(a_{x,u})<\phi_{\mu^x,u}(a_{x,u})$ where $a_{x,u}$ is
defined in~\eqref{E1}. Since $\phi_{\mu,u}\le
\phi_{\mu^x,u}\le\phi_{\nu, u}$, we obtain 
$\phi_{\mu,u}(a_{x,u})<\phi_{\nu,u}(a_{x,u})$, contradicting~\eqref{E1}.
\\ 

\noindent 
We proved that for every $u$, any point of the support of $K(x,dy)$
belongs to ${\cal D}_{x,u}$. This implies that the support of $K(x,dy)$ is
included in ${\cal D}_{x,u}$, achieving the proof.
\nprf
\begin{remark} 
 If \eqref{cnd01} holds for all $x\in \Supp (\mu )$ 
 and some $y_1,\ldots ,y_k \in \Supp (\nu)$ then we have 
\begin{equation} 
\label{condf} 
{\mathscr C} (\Supp (\mu ) ) \subset {\mathscr C} ( \Supp (\nu ) ), 
\end{equation} 
\end{remark} 
\bprf 
 Let $x \in \Supp (\mu )$. 
 If $\mu(x)\le \nu(x)$ then we clearly 
 have $x\in {\mathscr C} ( \Supp (\nu ))$, 
 and if $\mu(x) > \nu(x)$ then \eqref{cnd01} also implies 
 $x\in {\mathscr C} ( \Supp (\nu ))$. 
\nprf 
\noindent 
 As a consequence of the above remark 
 we note that the conclusion of Theorem~\ref{T5.1} is stronger than \eqref{cnd00}, 
 since two measures $\mu$ and $\nu$ may satisfy \eqref{condf} 
 without satisfying the condition $\mu \preceq_{\rm cxp} \nu$. 
 Counterexamples are easily constructed in dimension one when 
 $\mu$ and $\nu$ have same support. 
\\ 

\noindent 
 Let us now consider some examples with $d=2$ in the complex plane for 
 simplicity of notation. 
 First, let 
$$\mu=\frac12 (\delta_{-1}+\delta_1) \qquad 
 \mbox{and} 
 \qquad 
 \nu=\frac14(\delta_{-1-i}+\delta_{-1+i}+\delta_{1-i}+\delta_{1+i} ) 
, 
$$ 
 and $x=-1$. In this case we have $C_{-1}=[-1-i,-1+i]$, as illustrated in 
 the following figure: 
\\ 

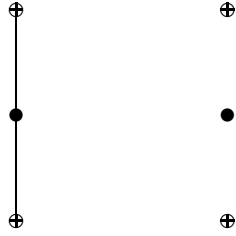
\begin{center} 
\begin{picture}(200,130)(-50,-100)
\linethickness{0.7pt}
\put(40,0){\circle*{5}} 
\put(-40,0){\circle*{5}} 
\put(-40,40){\circle{5}} 
\put(-42.5,40){\line(1,0){5}} 
\put(-40,42.5){\line(0,-1){5}} 
\put(-40,-40){\circle{5}} 
\put(-42.5,-40){\line(1,0){5}} 
\put(-40,-37.5){\line(0,-1){5}} 
\put(40,40){\circle{5}} 
\put(37.5,40){\line(1,0){5}} 
\put(40,42.5){\line(0,-1){5}} 
\put(40,-40){\circle{5}} 
\put(37.5,-40){\line(1,0){5}} 
\put(40,-37.5){\line(0,-1){5}} 
\linethickness{0.2pt}
\put(-40,-40){\line(0,1){80}} 
\end{picture}
\end{center} 
\vspace{-2cm} 
 In order to see this, take $u=1$ and then $u=-1$, 
 and note that $C_{-1}$ is contained in the vertical line 
 passing through $-1$ and apply Theorem~\ref{T5.1}. 
 Similarly we have $C_{1}=[1-i,1+i]$. 
 Next, consider 
$$ 
 \mu=\frac13\delta_{1/2}+\frac23\delta_{-1/4} 
 \qquad 
 \mbox{and} 
 \qquad 
 \nu=\frac13 (\delta_{1}+\delta_j+\delta_{j^2}) 
, 
$$ 
 where $j=e^{2i\pi/3}$. 
 Then Theorem~\ref{T5.1} shows that 
$$C_{1/ 2} = C_{-1/4} = 
 {\mathscr C} (\Supp (\nu )) 
, 
$$ 
 since $\{ y_1, y_2 , y_3\}$ is necessarily equal to $\Supp ( \nu )$, 
 as illustrated below: 
\\ 
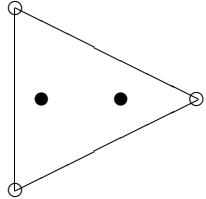
\begin{center} 
\begin{picture}(200,130)(-50,-100)
\linethickness{0.7pt}
\put(48.7,0){\circle{5}} 
\put(-20,34.5){\circle{5}} 
\put(-20,-34.2){\circle{5}} 
\put(20,0){\circle*{5}} 
\put(-10,0){\circle*{5}} 
\linethickness{0.2pt}
\put(-20,34.64){\line(2,-1){69.2}} 
\put(-20.1,34.64){\line(0,-1){69}} 
\put(-20,-34.64){\line(2,1){69.2}} 
\end{picture}
\end{center} 
\vspace{-2cm} 
\noindent Finally, note 
 that one can extend the result of Theorem~\ref{T5.1} to any non empty subset $E$ of $\real^d$ such that $x\in E$ implies 
$\mu(x)>\nu(x)$. 
 Letting 
$$ 
a_{E,u}=\inf\left\{b\in \left[\sup_{x\in E}\langle u,x\rangle,\infty\right[,\ \phi_{\mu,u}(b)=\phi_{\nu,u}(b)\right\} 
$$ 
 for $u\in S^{d-1}$, with $a_{E,u}=\sup_{x\in E}a_{x,u}$ since $E$ is finite, 
$$
{\cal D}_{E,u}= \{y\in \real^d \ : \ \langle u,y\rangle \le a_{E,u} \} 
, 
$$ 
 and defining 
$$
{\cal C}_E : =\bigcap_{u\in S^{d-1}}{\cal D}_{E,u} 
$$
 which is a compact convex set containing $E$, we have 
 ${\cal C}_E \subset \Supp (\nu)$ (since as previously $a_{E,u}\leq \tilde{a}_{u}$) and:  
\bcor
\label{C5.6} 
 Assume that $\mu$ and $\nu$ have finite supports and that 
 $\mu \preceq_{\rm cxp} \nu$. 
 Then for all non empty subset $E$ of $\real^d$ such that $x\in E$ implies 
$\mu(x)>\nu(x)$, there exists
 $k\in \{ 2, \ldots , (d+1){\rm card}(E)\}$ 
 and $k$ elements $\displaystyle y_1,\ldots,y_k \in \Supp (\nu)$ 
 distinct from $x$, such that 
$$ 
 E
 \ 
 \subset 
 \ 
 {\mathscr C} (\{y_1,\ldots,y_k\}) 
 \ 
 \subset 
 \ 
 {\cal C}_E 
. 
$$
\ncor
\bprf
From the definition of ${\cal C}_E$ it is clear that if $x\in E$ then ${\cal C}_x\subset {\cal C}_E$. Consequently applying Theorem~\ref{T5.1} to every $x\in E$ gives the result. 
\nprf
\noindent 
 When $\mu$ and $\nu$ are probability measures, 
 the existence of the admissible kernel $K$ such that $\mu K=\nu$, 
 used in the proof of Theorem~\ref{T5.1}, is also known as Strassen's theorem
 \cite{strassen}, and it is equivalent to the existence of 
 two random variables $F , G$ with respective laws
 $\mu$ and $\nu$, and such that $F = \ee [ G | F ]$. 
 Here we used Theorem~40 of \cite{delmey2} which relies on the Hahn-Banach theorem. 
 In dimension one this result has been recovered via a constructive proof 
 in \cite{muller}. 
 We close this paper with the following remark which concerns 
 the  $\preceq_{\rm cx}$ ordering. 
\begin{remark} 
 The conclusion of Theorem~\ref{T5.1}, associated to the 
 condition $\mu \preceq_{\rm cx} \nu$, 
 implies the existence of an admissible kernel $K$ such that 
 $\mu K=\nu$. 
\end{remark} 
\begin{Proof} 
 We use the notation of  Theorem~\ref{T5.1}. 
 First we show that if $\mu(\{x\})>\nu(\{x\})$ then 
 there exists a kernel $K_x$ such that $K_x(x,dy)$ is
supported by $\{x,y_1,\ldots,y_k\}$ and is not equal to $\delta_x$,
$K_x(x',dy)$ is equal to $\delta_{x'}$ if $x'\not=x$, 
and $\mu \preceq_{\rm cx} \mu K_x \preceq_{\rm cx} \nu$. 
 Indeed we can take 
$$
K_x(x,dy)=(1-\varepsilon)\delta_x(dy)+\varepsilon\sum_{i=1}^ka_i\delta_{y_i}(dy)
$$
where the $a_i$'s are positive, $\sum_{i=1}^ka_i=1$,
$x=\sum_{i=1}^ka_iy_i$ and $\varepsilon >0$ is sufficiently small 
 - the existence of $\varepsilon$ follows from the fact that the
functions $\phi_{\mu, u}$ and $\phi_{\nu, u}$ are continuous in $u$,
together with the compactness of $S^{d-1}$.
Now let $K$ be a maximal\footnote{Here, $K\le K'$ means 
that there exists an admissible  $K''$ such that $K\circ K''=K'$.} 
 admissible kernel such that $\mu K \preceq_{\rm cx} \nu$ and the support of $\mu K$ is included in $\Supp(\mu)\cup \Supp(\nu)$. If $\mu K\not=\nu$ then we can
apply the argument above to $\mu K$,  $\nu$ and $x$ such that $\mu
K(\{x\})>\nu(\{x\})$,  and find a non trivial
kernel $K_x$ such that $\mu KK_x \preceq_{\rm cx} \nu$,
contradicting the maximality of $K$. So we conclude that $\mu K=\nu$.
\end{Proof} 
\noindent 
 Thus an independent proof of Theorem~\ref{T5.1}, 
 not relying on Theorem~40 of \cite{delmey2}, would provide a direct 
 construction the admissible kernel $K$, extending the result 
 of \cite{muller} to higher dimensions. 

\small
\bibliographystyle{plain}
\def\cprime{$'$} \def\polhk#1{\setbox0=\hbox{#1}{\ooalign{\hidewidth
  \lower1.5ex\hbox{`}\hidewidth\crcr\unhbox0}}}
  \def\polhk#1{\setbox0=\hbox{#1}{\ooalign{\hidewidth
  \lower1.5ex\hbox{`}\hidewidth\crcr\unhbox0}}} \def\cprime{$'$}

\end{document}